\documentclass{article}

\usepackage{subfigure}
\usepackage{amsmath,amsfonts}
\usepackage{graphicx}

\def\d{\mathrm{d}}
\def\N{\mathrm{N}}
\def\R{\mathbb{R}}

\begin{document}

\title{Evolution Strategies in Optimization Problems\footnote{Partially presented at the
5th Junior European Meeting on \emph{Control and Information
Technology} (JEM'06), Sept 20--22, 2006, Tallinn, Estonia. Research Report CM06/I-44.}}

\author{Pedro A. F. Cruz \and Delfim F. M. Torres}

\date{Department of Mathematics\\
      University of Aveiro\\
      3810-193 Aveiro, Portugal\\
      \medskip
      \texttt{\{pedrocruz,delfim\}@ua.pt}}

\maketitle


\begin{abstract}
Evolution Strategies are inspired in biology and part of a larger
research field known as Evolutionary Algorithms. Those strategies
perform a random search in the space of admissible functions,
aiming to optimize some given objective function. We show that
simple evolution strategies are a useful tool in optimal control,
permitting to obtain, in an efficient way, good approximations to
the solutions of some recent and challenging  optimal control
problems.
\end{abstract}

\smallskip


\smallskip


\smallskip


\medskip


\section{Introduction}

Evolution Strategies (ES) are algorithms inspired in biology, with
publications dating back to $1965$ by separate authors
H.P.~Schwefel and I.~Rechenberg (\textrm{cf.} \cite{bib:ev}). ES
are part of a larger area called Evolutionary Algorithms that
perform a random search in the space of solutions aiming to
optimize some objective function. It is common to use biological
terms to describe these algorithms. Here we make use of a simple
ES algorithm known as the $(\mu,\lambda)-$ES method \cite{bib:ev},
where $\mu$ is the number of \textit{progenitors} and $\lambda$ is
the number of generated approximations, called
\textit{offsprings}. Progenitors are \textit{recombined} and
\textit{mutated} to produce, at each generation, $\lambda$
\textit{offsprings} with innovations sampled from a multivariate
normal distribution. The variance can also be subject to mutation,
meaning that it is part of the \textit{genetic code} of the
population. Every solution is evaluated by the objective function
and one or some of them selected to be the next
\textit{progenitors}, allowing the search to go on, stopping when
some criteria is met. In this paper we use a recent convergence
result proved by A.~Auger in 2005 \cite{bib:auger}. The log-linear
convergence is achieved for the optimization problems
we investigate here, and depends on the number $\lambda$ of search points.

Usually optimal control problems are approximately solved by means of numerical
techniques based on the gradient vector or the Hessian matrix
\cite{Viorel}. Compared with these techniques, ES provide easier
computer coding because they only use measures from a discretized
objective function. A first work combining these two research
fields (ES and optimal control) was done by D.S.~Szarkowicz in 1995
\cite{bib:szar}, where a Monte Carlo method (an algorithm with the
same principle as ES) is used to find an approximation to the
classical brachistochrone problem. In the late nineties of the XX
century, B.~Porter and his collaborators showed how ES are useful
to synthesize optimal control policies that minimize manufacturing
costs while meeting production schedules
\cite{ManufacturingSystems}. The use of ES in Control has grown
during the last ten years, and is today an active and promising
research area. Recent results, showing the power of ES in Control,
include Hamiltonian synthesis \cite{Kloster2001}, robust
stabilization \cite{Paris2000}, and optimization
\cite{OptimalElevator}. Very recently, it has also been shown that
the theory of optimal control provides insights that permit to
develop more efective ES algorithms \cite{OChelpsES}.

In this work we are interested in two classical problems of the
calculus of variations: the 1696 brachistochrone problem and the
1687 Newton's aerodynamical problem of minimal resistance (see
\textrm{e.g.} \cite{Tikh}). These two problems, although
classical, are source of strong current research on optimal
control and provide many interesting and challenging open
questions \cite{bib:plakhov,brachisSussmann}. We focus our study
on the brachistochrone problem with restrictions proposed by
A.G.~Ramm in 1999 \cite{bib:ramm}, for which some questions still
remain open (see some conjectures in \cite{bib:ramm}); and on a
generalized aerodynamical minimum resistance problem with
non-parallel flux of particles, recently studied by Plakhov and
Torres \cite{bib:plakhov,bib:delfim}. Our results show the
effectiveness of ES algorithms for this class of problems and
motivate further work in this direction in order to find the (yet)
unknown solutions to some related problems, as the ones formulated
in \cite{Billiards}.


\section{Problems and Solutions}

All the problems we are interested in share the same formulation:
  \[
  \min T[y(\cdot)] = \int_{x_0}^{x_f} L(x,y(x),y'(x)) \d x
  \]
on some specified class of functions, where $y(\cdot)$ must
satisfy some given boundary conditions $(x_0,y_0)$ and
$(x_f,y_f)$.

We consider a simplified $(\mu,\lambda)-$ES algorithm where we put
$\mu=1$, meaning that on each generation we keep only one
\textit{progenitor} to generate other candidate solutions, and set
$\lambda=10$ meaning we generate $10$ candidate solutions called
\textit{offsprings} (this value appear as a reference value in the
literature). Also, the algorithm uses an individual and constant
$\sigma^2$ variance on each coordinate, which is fixed to a small
value related with the desired precision. The number of iterations
was $100\,000$ and $\sigma^2$ was tuned for each problem. We got
convergence in useful time. The simplified $(1,10)-$ES algorithm
goes as follow:
\begin{enumerate}

\item Set an equal spaced sequence of $n$ points
$\{x_0,\ldots,x_i,\ldots,x_f\}$ where $i=1,\ldots,n-2$; $x_0$ and
$x_f$ are kept fixed (given boundary conditions);

\item Generate a randomly piecewise linear function $y(\cdot)$
that approximate the solution, defined by a vector $y =
\{y_0,\ldots,y_i,\ldots,y_f\},\,i=1,\ldots,n-2$; transform $y$ in
order to satisfy the boundary conditions $y_0$ and $y_f$ and the
specific problem restrictions on $y$, $y'$ or $y''$;

\item Do the following steps a fixed number $N$ of times:
\begin{enumerate}
  \item based on $y$ find $\lambda$ new candidate solutions
  $Y^c$, $c=1,\ldots,\lambda$, where each new candidate is
  produced by $Y^c = y + \N(0,\sigma^2)$ where $\N(0,\sigma^2)$
  is a vector of random perturbations
  from a normal distribution; transform each $Y^c$ to obey
  boundary conditions $y_0$ and $y_f$ and
  other problem restrictions on $y$, $y'$ or $y''$;
  \item determine $T^c := T[Y^c]$, $c=1,\ldots,\lambda$,
  and choose the new $y:=Y^c$ as the one with minimum $T^c$.
\end{enumerate}

\end{enumerate}
In each iteration the best solution must be kept because
$(\mu,\lambda)-$ES algorithms don't keep the best solution from
iteration to iteration.

The next subsections contain a description of the studied
problems, respective solutions and the approximations found by the
described algorithm.

\subsection{The classical brachistochrone problem, 1696}
\label{subsec:ClassBrach}


\textbf{Problem statement.}~The brachistochrone problem consists
in determining the curve of minimum time when a particle starting
at a point $A=(x_0,y_0)$ of a vertical plan goes to a point
$B=(x_1,y_1)$ in the same plane under the action of the gravity
force and with no initial velocity. According to the energy
conservation law $\frac{1}{2} m v^2 + m g y = m g y_0$ one easily
deduce that the time a particle needs to reach $B$ starting from
point $A$ along curve $y(\cdot)$ is given by
  \begin{equation}\label{eq:bra}
  T[y(\cdot)] = \frac{1}{\sqrt{2g}} \int_{x_0}^{x_1} \sqrt{\frac{1 + (y')^2}{y_0 - y}} \d x
  \end{equation}
where $y(x_0)=y_0$, $y(x_1)=y_1$, and $y \in C^2(x_0,x_1)$. The
minimum to \eqref{eq:bra} is given by the famous Cycloid:
  \begin{equation*}
  \gamma\,:\,\left\{
    \begin{array}{l}
    x = x_0 + \frac{a}{2}(\theta - \sin \theta) \\
    y = y_0 - \frac{a}{2}(1 - \cos \theta)
    \end{array}
  \right.
  \end{equation*}
with $\theta_0 \le \theta \le \theta_1$, $\theta_0$ and $\theta_1$
the values of $\theta$ in the starting and ending points
$(x_0,y_0)$ and $(x_1,y_1)$. The minimum time is given by
$T=\sqrt{a/(2g)}\theta_1$, where parameters $a$ and $\theta_1$ can
be determined numerically from boundary conditions.

\begin{figure}
\centering \subfigure[\label{fig:bra-a}Continuous line with dots
is the piecewise approximate solution; the dashed line the optimal
solution.]{
\includegraphics[width=5.5cm,height=6cm]{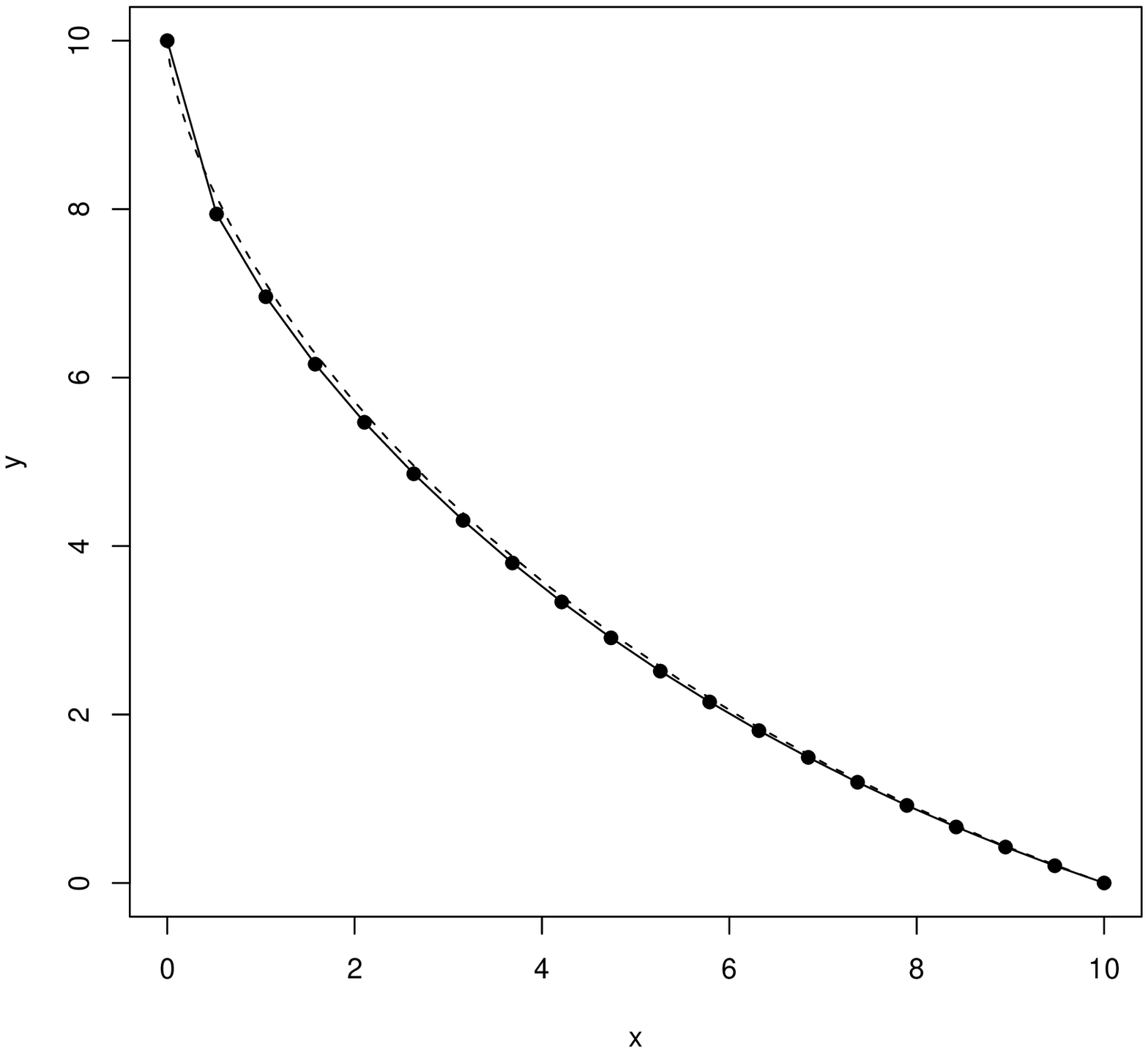}
}~ \subfigure[\label{fig:bra-b}Logarithm of iterations vs.
logarithmic distance to the minimum value of \eqref{eq:bra}]{
\includegraphics[width=5.5cm,height=6cm]{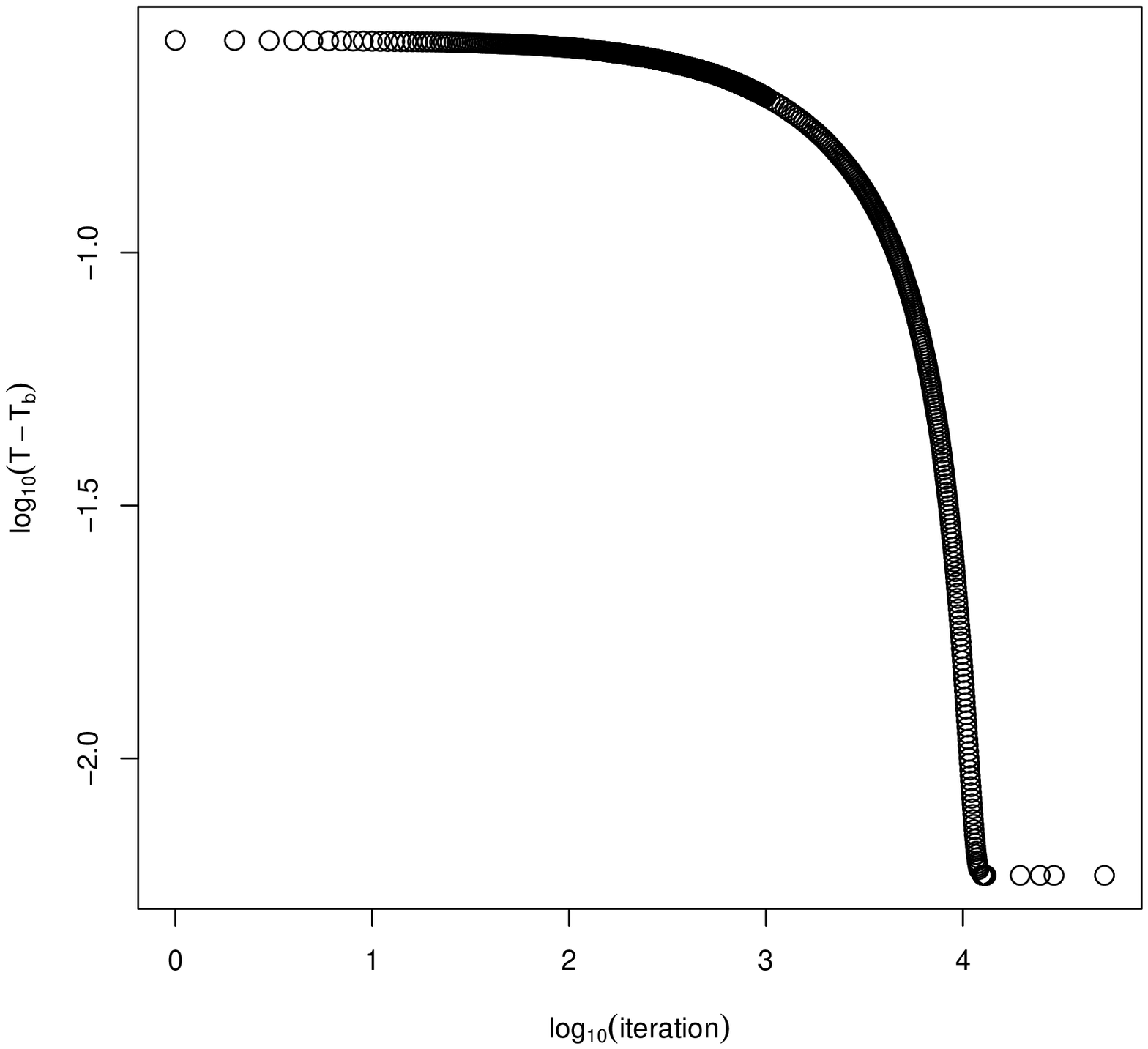}
}\\
\caption{The brachistochrone problem and approximate solution.}\label{fig:br}
\end{figure}

\textbf{Results and implementation details.}~Consider the
following three curves and the correspondent time a particle needs
to go from $A$ to $B$ through them:
\begin{description}
\item[$T_b$:]~The brachistochrone for the problem with
$(x_0,x_1)=(0,10)$, $(y_0,y_1)=(10,0)$ has parameters $a\simeq
5.72917$ and $\theta_1 =2.41201$; the time is $T_b\simeq1.84421$;
\item[$T_{es}$:]~A piecewise linear function with $20$ segments shown
in fig.~\ref{fig:bra-a} was found by ES; the time is $T_{es}=1.85013$;
\item[$T_o$:]~A piecewise linear function with $20$ segments
defined over the Brachisto\-chrone; the time is $T_o=1.85075$.
\end{description}
From fig.~\ref{fig:bra-a} one can see that the piecewise linear solution
is made of points that are not over the brachistochrone because
that is not the best solution for piecewise functions.
We use $\sigma=0.01$ (see appendix for cpu-times).
Fig.~\ref{fig:bra-b} shows that a little more than $10\,000$
iterations are needed to reach a good solution for the $20$ line
segment problem.

\subsection{Brachistochrone problem with restrictions, 1999}
\label{subsec:RammBrach}

\textbf{Problem statement.}~Ramm (1999) \cite{bib:ramm} presents a
conjecture about a brachistochrone problem over the set $S$ of
convex functions $y$ (with $y''(x)\ge0$ a.e.) and $0 \le y(x) \le
y_0(x)$, where $y_0$ is a straight line between $A=(0,1)$ and
$B=(b,0)$, $b>0$. Up to a constant, the functional to be minimized
is formulated as in (\ref{eq:bra}):
  \begin{equation*}
  T[y(\cdot)] = \int_0^b \frac{ \sqrt{1+(y')^2} }{ \sqrt{1-y} }\d
  x\, .
  \end{equation*}
Let $P$ be the line connecting $AO$ and $OB$, where
$O=(0,0)$; $P_{br}$ be the polygonal line connecting $AC$ and
$CB$, $C=(\pi/2,0)$. Then, $T_0:=T(y_0)=2\sqrt{1+b^2}$,
$T_P:=T(P)=2+b$, $T(P_{br})=\sqrt{4+\pi^2}+b-\pi/2$.
Let the brachistochrone be $y_{br}$. The following inequalities,
for each $y \in S$, hold \cite{bib:ramm}:
\begin{enumerate}
\item if $0 < b < 4/3$ then $T(y_{br}) \le T(y) < T_P$;
\item if $4/3 \le b \le \pi/2$ then $T(y_{br}) \le T(y) \le T_0$;
\item if $b > \pi/2$ then $T(P_{br}) < T(y) \le T_0$.
\end{enumerate}

The classical brachistochrone solution holds for cases 1 and 2
only. For the third case, Ramm has conjectured that the minimum
time curve is composed by the brachistochrone between $(0,1)$ and
$(\pi/2,0)$ and then by the horizontal segment between $(\pi/2,0)$
and $(x_f,0)$.

\begin{figure}
\centering \subfigure[\label{fig:ramm-a}Continuous line with dots
is the obtained approximated solution; dashed line Ramm's
conjectured solution.]{
          \includegraphics[width=5.5cm,height=6cm]{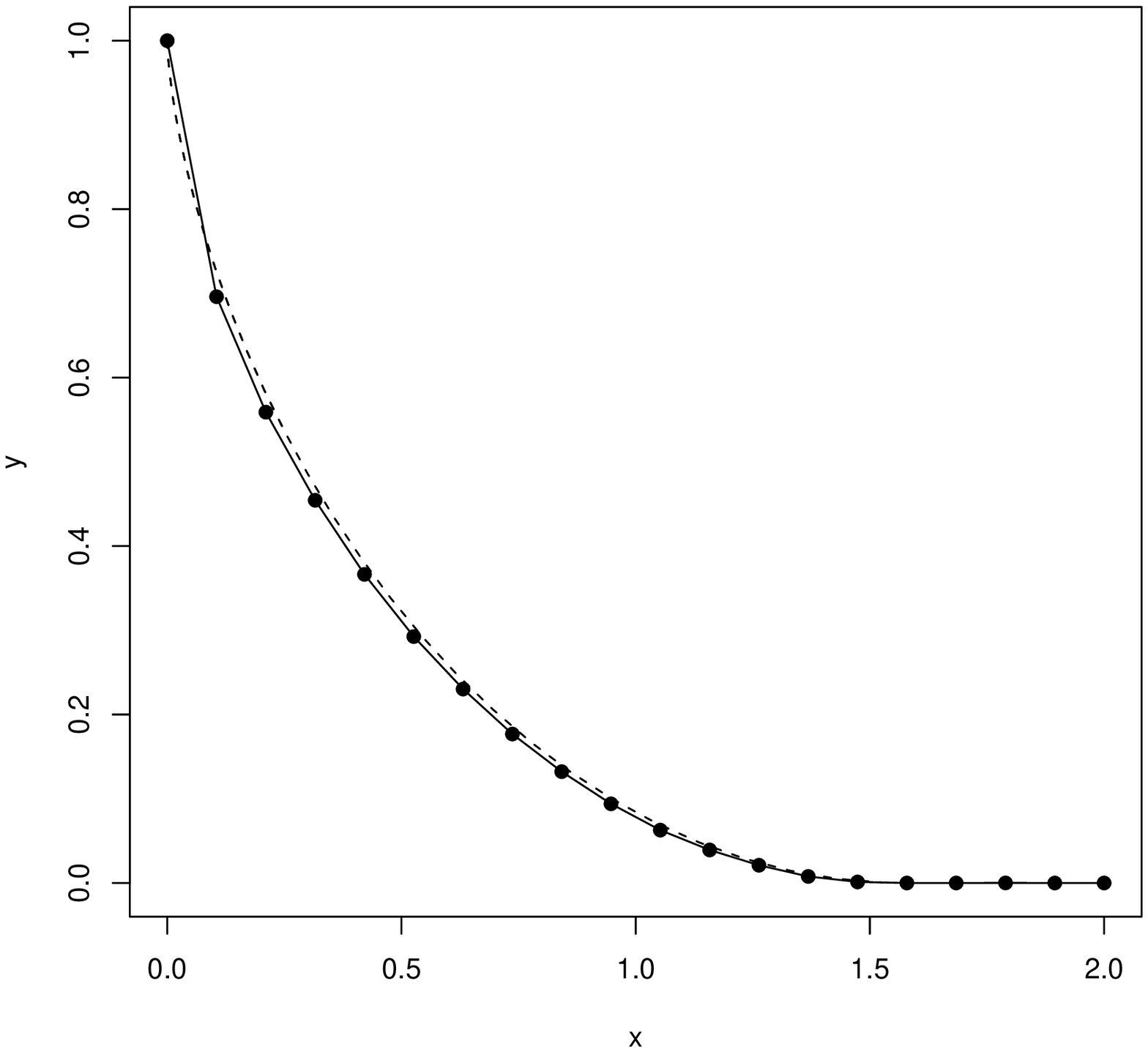}}~
\subfigure[\label{fig:ramm-b}Logarithm of iterations vs.
logarithmic distance to minimum integral value.]{
          \includegraphics[width=5.5cm,height=6cm]{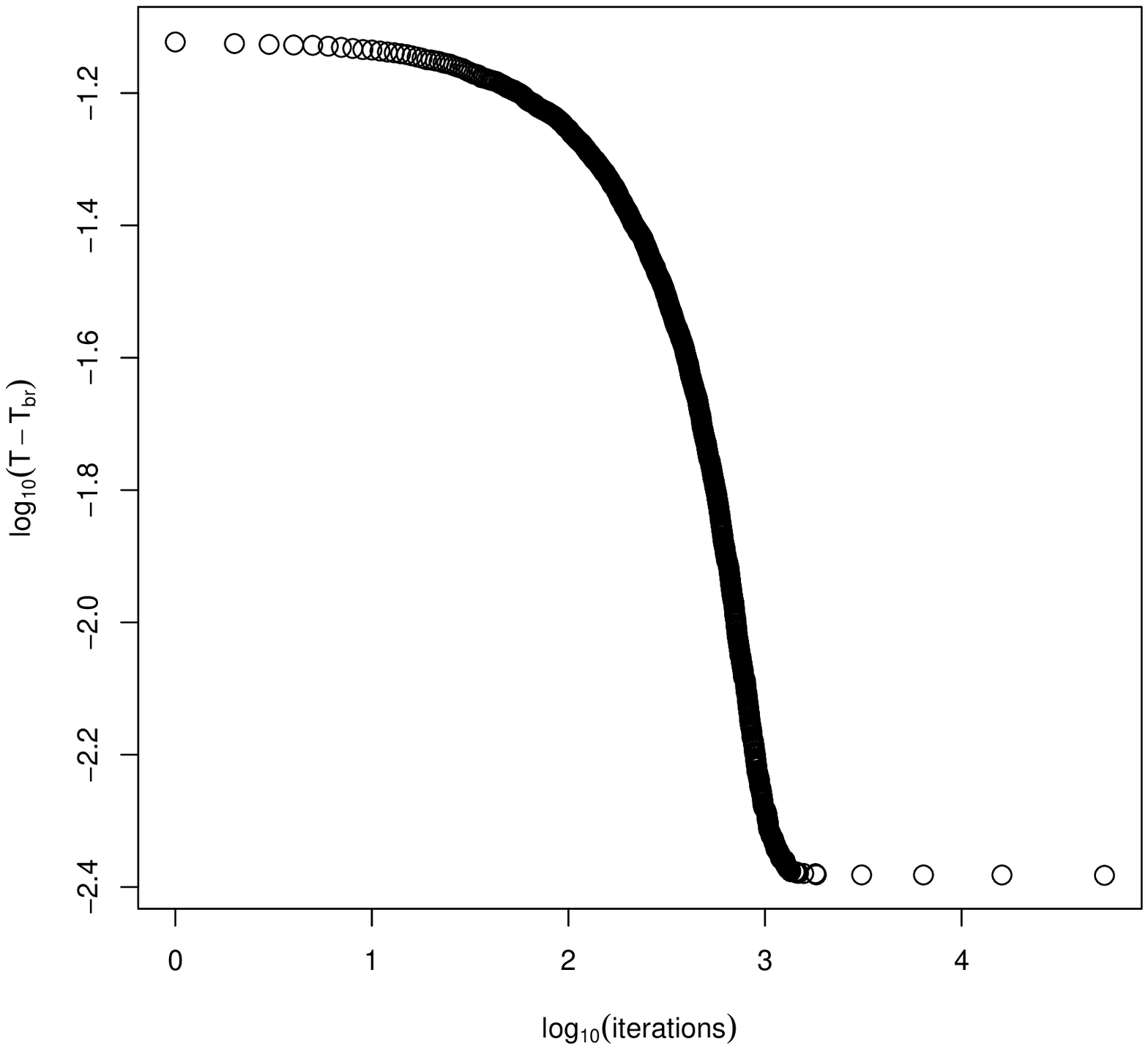}
}\\
\caption{Ramm's conjectured solution and approximate solution.}\label{fig:ramm}
\end{figure}

\textbf{Results and implementation details.}~We study the problem
with $b=2$. Our results give force to Ramm's conjecture mentioned
above for case 3. We compare three descendant times:
\begin{description}
\item[$T_{br}$:]~The conjectured solution in continuous time
takes $T_{br} = \sqrt{\alpha/9.8} \theta_f + (b-\pi/2) /
\sqrt{2*9.8} = 0.8066$; \item[$T_{es}$:]~The $20$ segment
piecewise linear solution found by ES needs $T_{es}=0.8107$;
\item[$T_o$:]~The $20$ segment piecewise linear solution with
points over the conjectured solution needs $T_o=0.8111$.
\end{description}
Previous values and fig.~\ref{fig:ramm} permits to take similar
conclusions than the ones obtained for the pure brachistochrone
problem (\S\ref{subsec:ClassBrach}).
We use $\sigma=0.001$ (see appendix for cpu-times).
Fig.~\ref{fig:ramm-b} shows that less than $10\,000$ iterations
are needed to reach a good solution.

\subsection{Newton's minimum resistance, 1687}
\label{subsec:NewtonClassical}

\textbf{Problem statement.}~Newton's aerodynamical problem
consists in determining the minimum resistance profile of a body
of revolution moving at constant speed in a rare medium of equally
spaced particles that don't interact with each other. Collisions
with the body are assumed to be perfectly elastic. Formulation of
this problem is: minimize
  \begin{equation*}
  R[y(\cdot)] = \int_0^r
  \frac{ x }{ 1 + \dot y(x)^2 }  \d x
  \end{equation*}
where $0\le x \le r$, $y(0)=0$, $y(r)=H$ and $y'(x)\ge0$.
The solution is given in parametric form:
  \begin{equation*}
  x(u)  =  2\lambda u \, , \quad
  y(u)  =  0 \, , \quad \text{ for } u \in[0,1] \, ;
  \end{equation*}
  \begin{equation*}
  x(u) =  \frac{\lambda}{2}(\frac{1}{u} + 2u + u^3) \, , \quad
  y(u) =  \frac{\lambda}{2}(-\log u + u^2 + \frac{3}{4} u^4) - \frac{7
  \lambda}{8}\, , \quad \text{ for } u \in [1, u_\textrm{max}] \, .
  \end{equation*}
Parameters $\lambda$ and $u_\textrm{max}$ are obtained solving
$x(u_{\max})=r$ and $y(u_{\max})=H$.

\bigskip\bigskip

\begin{figure}
\centering \subfigure[\label{fig:newton-a}Continuous line with dots
is the obtained approximation; the dashed line the optimal solution.]{
          \includegraphics[width=5.5cm,height=6cm]{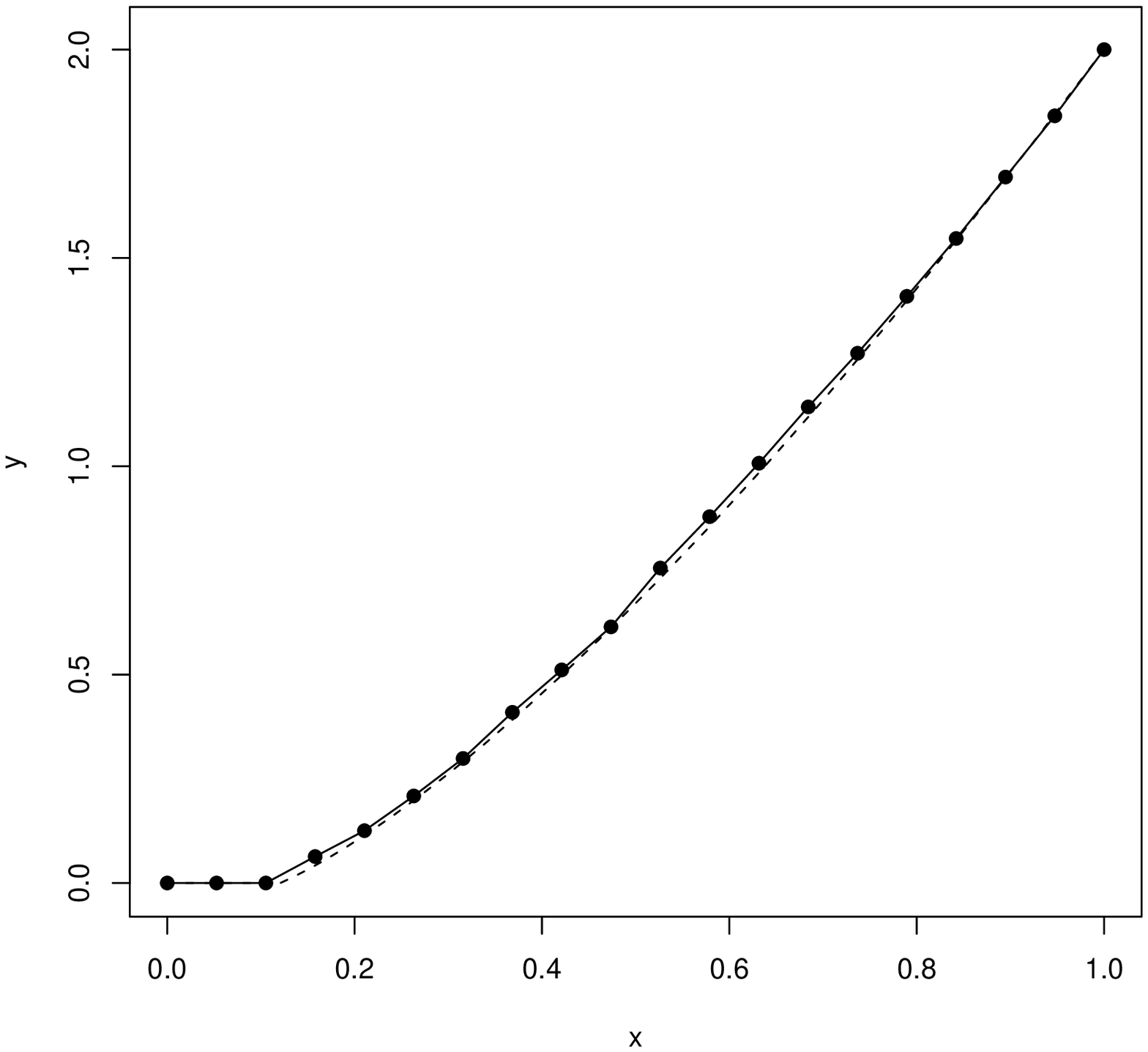}
}~ \subfigure[\label{fig:newton-b}Logarithm of iterations vs.
logarithmic distance to minimum integral value.]{
          \includegraphics[width=5.5cm,height=6cm]{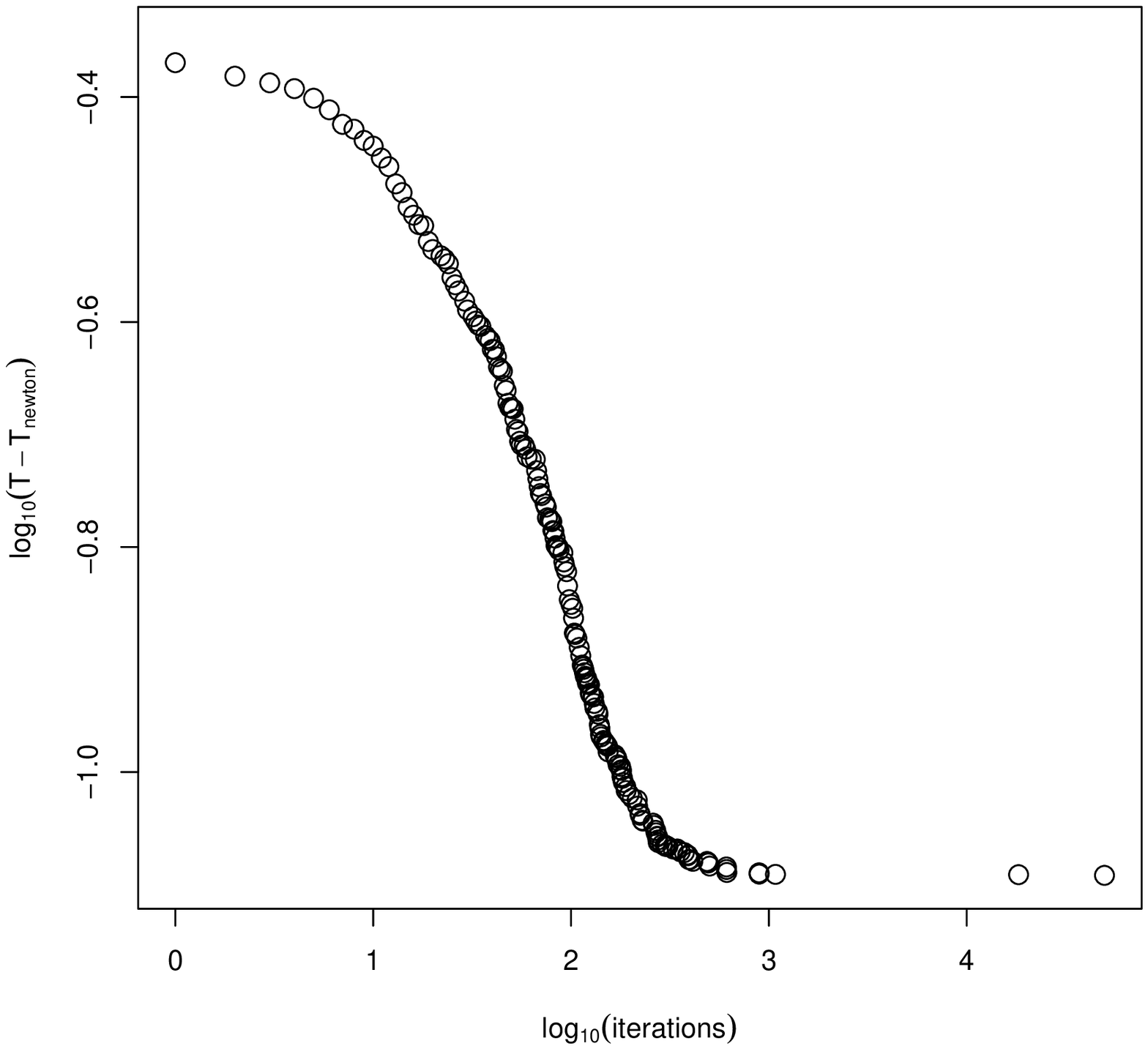}
}\\
\caption{Optimal solution to Newton's problem and
approximation.}\label{fig:newton}
\end{figure}

\textbf{Results and implementations details.}~For $H=2$ we have:
\begin{description}
\item[$R_{newton}$:]~The exact solution has resistance $R_{newton} = 0.0802$;
\item[$R_{es}$:]~The $20$ segment piecewise linear solution found by ES has $R_{es}=0.0809$;
\item[$R_o$:]~The $20$ segment piecewise linear solution with points over the exact
solution leads to $R_o=0.0808$.
\end{description}
Newton's problem reveals to be more complex than previously
studied brachistochrone problems. Trial-and-error
was needed in order to find a useful $\sigma^2$ value.
For example, using $\sigma=0.001$ our algorithm
seems to stop in some local minimum. In
fig.~\ref{fig:newton} an approximate solution with $\sigma=0.01$ is shown.
We also have observed that changing the starting point causes minor differences in the
approximate solution. The achieved ES solution should be better since $R_o$ is
better than $R_{es}$. One possible explanation for this fact
is that we are using $20$ $x_i$ fixed points and
the optimal solution has a break point at $x=2\lambda$.
We use $\sigma=0.01$ (see appendix for cpu-times).
Fig.~\ref{fig:newton-b} shows that less than $1\,000$
iterations are needed to reach a good solution.

\subsection{Newton's problem with temperature, 2005}
\label{subsec:NewtonTemp}

\textbf{Problem statement.}~The problem consists in determining
the body of minimum resistance, moving with constant velocity in a
rarefied medium of chaotically moving particles with velocity
distributions assumed to be radially symmetric in the Euclidian
space $\mathbb{R}^d$. This problem was posed and solved in
2005-2006 by Plakhov and Torres \cite{bib:plakhov,bib:delfim}. It
turns out that the two-dimensional problem ($d=2$) is more richer
than the three-dimensional one, being possible five types of
solutions when the velocity of the moving body is not 'too slow'
or 'too fast' compared with the velocity of particles.

The pressure at the body surface is described by two functions: in
the front of the body the flux of particles causes resistance, in
the back the flux causes acceleration. We consider functions found
in \cite{bib:delfim}, where the two flux functions $p_+$ and $p_-$
are given by $p_+(u)=\frac{1}{1+u^2}+0.5$ and $p_-(u) =
\frac{0.5}{1+u^2}-0.5$. We also consider a body of fixed radius
$1$. The optimal solution depends on the body height $h$: the
front solution is denoted by $f_{h_+}$, which depends on some
appropriate front height $h_+$; and the solution for the rear is
denoted by $f_{h_-}$, depending on some appropriate height $h_-$.
Optimal solutions $f_{h_+}$ and $f_{h_-}$ are obtained:
  \begin{equation*}
  f_{h_+} = \min_{f_h} \quad \R_+(f_h) = \int_0^1 p_+(f_h'(t)) \d t
  \end{equation*}
and
  \begin{equation*}
  f_{h_-} = \min_{f_h}  \R_-(f_h) = \int_0^1 p_-(f_h'(t)) \d t.
  \end{equation*}
Then, the body shape is determined by minimizing
  \begin{equation*}
  R(h) = \min_{h_+ + h_-=h} ( \R_+(f'_{h_+}) + \R_-(f'_{h_-}) ) \,
  .
  \end{equation*}
Solution can be of five types ($d = 2$). From functions $p_+$ and
$p_-$ one can determine constants $u_+^0$, $u_*$, $u_-^0$ and
$h_-$. Then, depending on the choice of the height $h$, theory
developed in \cite{bib:plakhov,bib:delfim} asserts that the
minimum resistance body is:
  \begin{enumerate}
  \item a trapezium if $0 < h < u_+^0$;
  \item an isosceles triangle if $u_+^0 \le h \le u_*$;
  \item the union of a triangle and a trapezium if $u_* < h < u_* + u_-^0$;
  \item if $h \ge u_* + u_-^0$ the solution depends on $h_-$ and can be a union of
        two isosceles triangles with common base with heights $h_+$ and $h_-$ or
        the union of two isosceles triangles and a trapezium;
  \item a combination of a triangle, trapezium and other triangle, depending
        on some other particular conditions (\textrm{cf.}
        \cite{bib:plakhov}).
  \end{enumerate}

\textbf{Results and implementation details.}~We illustrate the
use of ES algorithms for $h=2$. Following section
4.1 of \cite{bib:plakhov} we have $u_*\simeq1.60847$ and $u_-^0=1$,
so this is case 3 above: $u_* < h < u_* + u_-^0$. The resistance values
are:
\begin{description}
\item[$R_{pd}$:]~The exact solution has resistance $R_{pd} = 0.681$;
\item[$R_{es}$:]~The $31$ segment piecewise linear
solution found by ES has $R_{es}=0.685$;
\end{description}
Similar to the classical problem of Newton (\S\ref{subsec:NewtonClassical}),
some hand search for the parameter $\sigma^2$ was needed.
\begin{figure}
\centering \subfigure[\label{fig:plakhov-a}Continuous line with
dots is the obtained approximation; dashed line the optimum.]{
          \includegraphics[width=5.5cm,height=6cm]{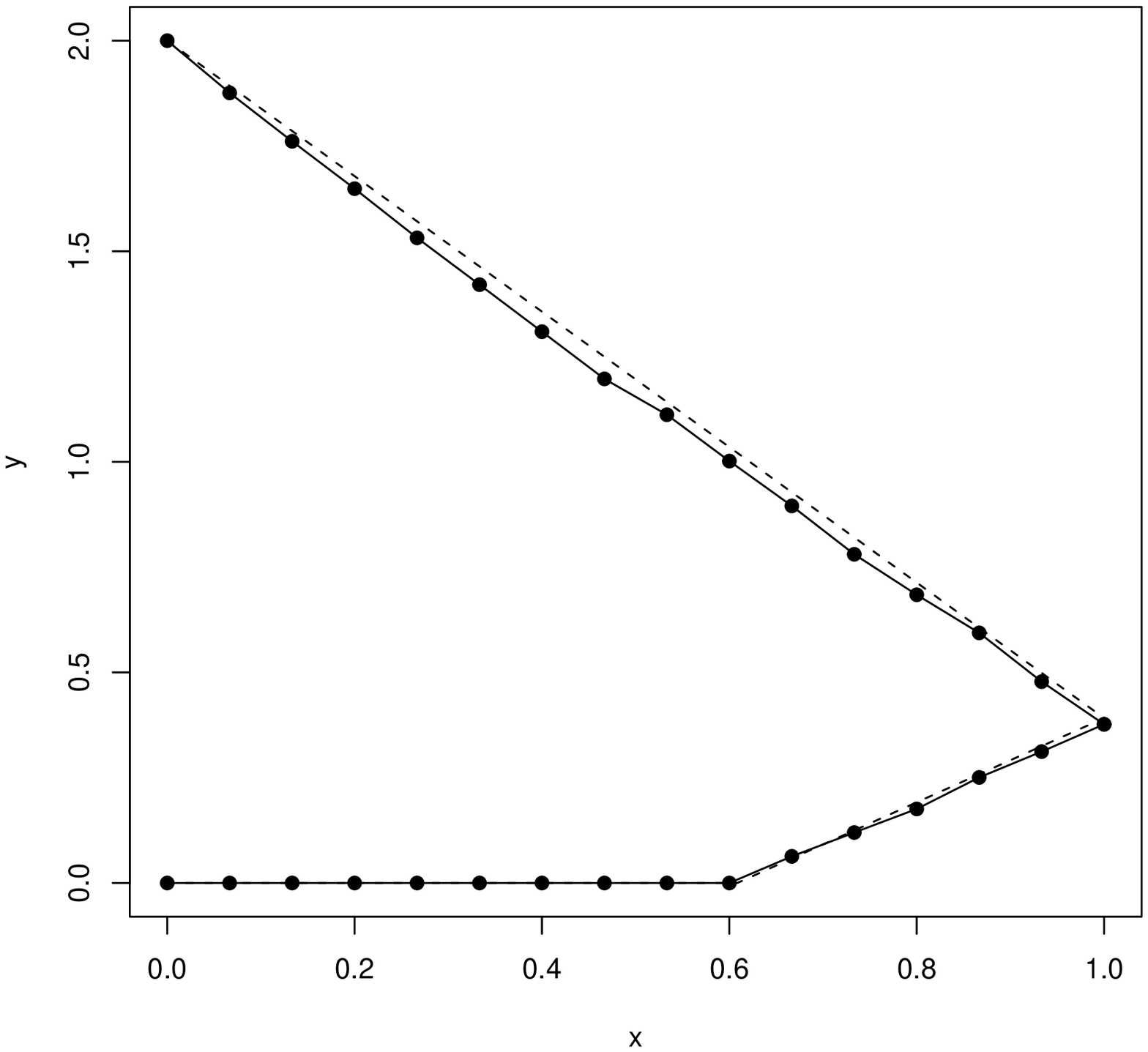}
}~ \subfigure[\label{fig:plakhov-b}Logarithm of iterations vs.
logarithmic distance to minimum integral value.]{
          \includegraphics[width=5.5cm,height=6cm]{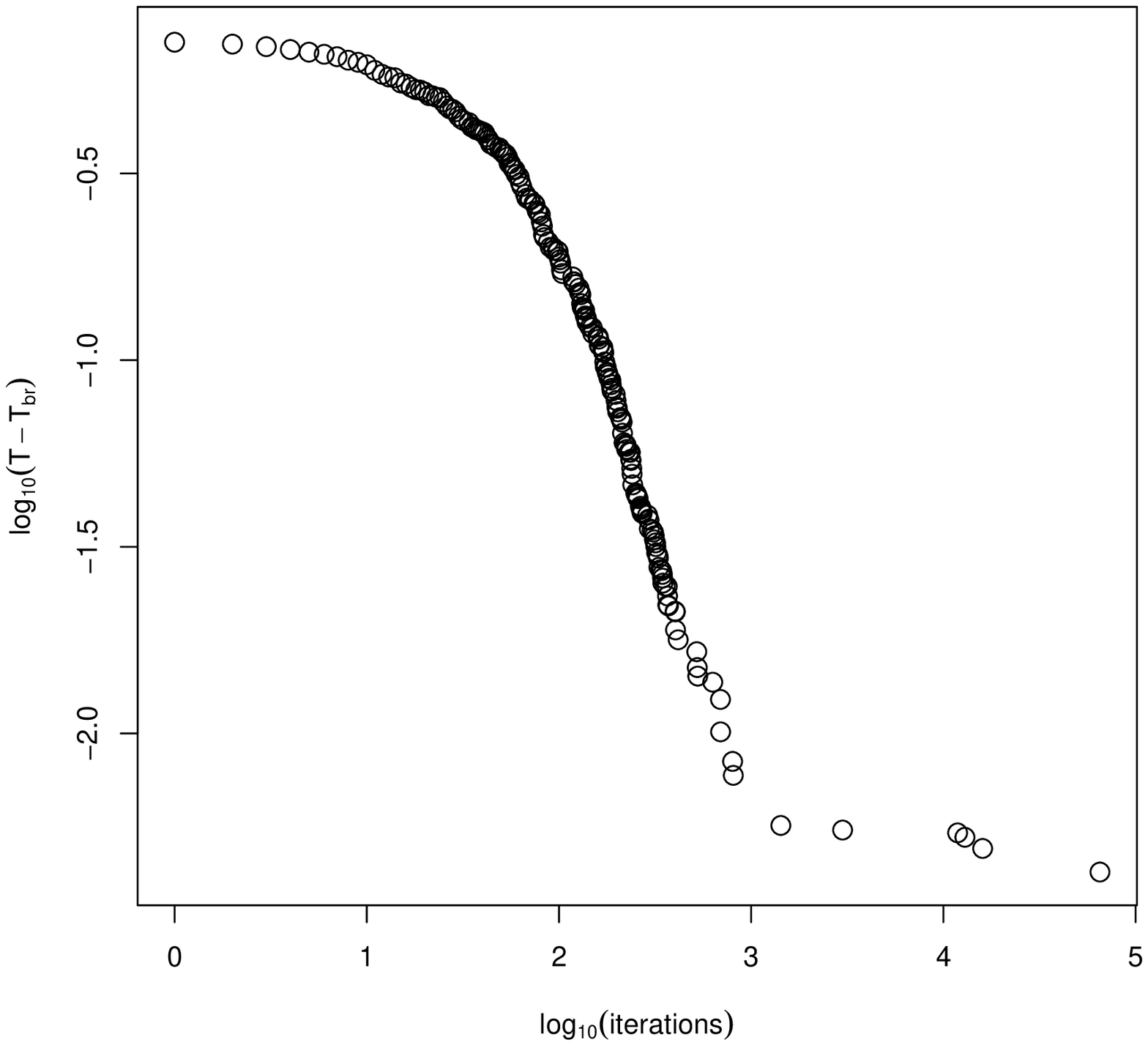}
}\\
\caption{2D Newton-type problem with
temperature.}\label{fig:plakhov}
\end{figure}
We use $\sigma=0.01$ and piecewise approximation with $31$ equal
spaced segments in $xx$ (see appendix for cpu-times).
Fig.~\ref{fig:newton-b} shows that only little more than $1\,000$
iterations are needed to reach a good solution.


\section{Conclusions and future directions}

Our main conclusion is that a simple ES algorithm can be
effectively used as a tool to find approximate solutions to some
optimization problems. In the present work we report
simulations that motivate the use of ES algorithms to find good
approximate solutions to brachistochrone-type and Newton-type
problems. We illustrate our approach with the classical problems
and with some recent and still challenging problems. More
precisely, we considered the 1696 brachistochrone problem (B); the
1687 Newton's aerodynamical problem of minimal resistance (N); a
recent brachistochrone problem with restrictions (R) studied by
Ramm in 1999, and where some open questions still remain
\cite{bib:ramm}; and finally a generalized aerodynamical minimum
resistance problem with non-parallel flux of particles (P),
recently studied by Plakhov and Torres
\cite{bib:plakhov,bib:delfim} and which gives rise to other
interesting questions \cite{Billiards}.

We argue that the approximated solutions we have found by the ES
algorithm are of good quality. We give two reasons. First, for the
Brachistochrone and Ramm's problems the functional value for the
ES approximation was better than the linear interpolation over the
exact solution, showing that ES algorithm is capable of a good
precision. The second reason is the low relative error
$r(T_Y,T_y)$ between the functional over the exact solution $T_y$
and the approximate solution $T_Y$, as shown in the following
table:
\begin{center}
\begin{tabular}{|l|c|c||l|c|c|}
\hline
Pr. & $\max|Y_k - y_k|$ & $r(T_Y,T_y)$ & Pr.& $\max|Y_k - y_k|$ & $r(T_Y,T_y)$ \\
\hline
(B) & 0.15 & 0.001 &
(N) & 0.08 & 0.01 \\
(R) & 0.09 & 0.003 &
(P) & 0.07 & 0.001 \\
\hline
\end{tabular}
\end{center}
where $y_k$ are points over the exact solution of the problem and
$Y_k$ are points from the piecewise approximation. We note that
$\max|Y_k - y_k|$ need not to be zero because the best continuous
solution and the best linear solution cannot be superposed.

ES algorithms use computers in an intensive way. For
brachistochrone-type and Newton-type problems, and nowadays
computing power, few minutes of simulation (or less) were enough
on an interpreted language (see appendix).

More research is needed to tune this kind of algorithms and obtain
more accurate solutions. Special attention must be put in
qualifying an obtained ES approximation: Is it a minimum of the
energy function? Is it local or global? Another question is
computer efficiency. Waiting few minutes in recent computers is
not bad, but can we improve the running times?

Concerning the accuracy, several new ES algorithms have been
proposed. These algorithms can tune $\sigma$ values and use
generated second order information that can influence the
precision and time needed. Also the use of random $xx$ points
(besides $y$ piecewise linear solution) should be investigated.

We believe that the simplicity of the technique considered in the
present work can help in the search of solutions to some open
problems in optimal control. This is under investigation and will
be addressed elsewhere.


\section*{Appendix -- hardware and software}

The code developed for this work can be freely obtained from the
first author's web page, at
\verb+http://www.mat.ua.pt/jpedro/evolution/+.

In most of our investigations few minutes were sufficient for
getting a good approximation for all the considered problems, even
using a code style prone to humans rather than machines (code was
done concerning clearness of concepts rather than execution
speed). Our simulations used a Pentium~4 CPU 3~GHz, running Debian
Linux \verb+http://www.debian.org+. The language was R
\cite{bib:r}, chosen because it is a fast interpreted language,
numerically oriented to statistics and freely available.

The following CPU-times were obtained with command
$$
\verb+time R CMD BATCH problem.R+
$$
where \verb+time+ keeps track of cpu used and $\verb+R+$ calls the
interpreter. The times are rounded and the last column estimates
the time for a first good solution:
\begin{center}
\begin{tabular}{lccc}
Problem & Section & $100\,000$ iterations & `Good solution' at \\
\hline
Brachistochrone & \S\ref{subsec:ClassBrach} & 10~min & 1~min\\
Ramm conjecture & \S\ref{subsec:RammBrach} & 10~min & 1~min\\
Newton & \S\ref{subsec:NewtonClassical} & 9~min & 10~sec\\
Plakhov \& Torres & \S\ref{subsec:NewtonTemp} &  14~min & 10~sec\\
\hline
\end{tabular}
\end{center}
We note that the per iteration `step' was $\sigma=0.001$ in the
brachistochrone(-type) problems and $\sigma=0.01$ for the
Newton(-type) problems. Using a compiled language like \textsf{C}
one can certainly improve times by several orders of magnitude.


\section*{Acknowledgements}

The authors were supported by the R\&D unit {\it Centre for
Research on Optimization and Control} (CEOC) from the Portuguese
Foundation for Science and Technology, cofinanced by the European
Community fund FEDER/POCI 2010.



\end{document}